# Multiple Object Trajectography

## Using Particle Swarm Optimization Combined to Hungarian Method

Max CERF [*]

Ariane Group, France

**Abstract**

The problem of simultaneous trajectography of several dynamical objects is formulated as an optimization problem. The available observations consist in a series of photographs showing undiscriminated objects. The goal is to find the object initial states so that the resulting trajectories match as well as possible the set of observations. An assignment problem is solved at each observation date by the Hungarian method, yielding a deviation cost between the simulated trajectories and the measurements. A fitness function summing the deviation costs is minimized by a particle swarm algorithm. The method is illustrated on a space orbitography application.

## 1. Introduction

Trajectography of dynamical objects is a common task in various scientific and engineering disciplines. The goal is to find initial conditions and model parameters so that the resulting simulated trajectory match a set of measurements. The estimated model can then be used for predictions and tracking purposes. References [1,2] provide a comprehensive presentation of standard methods like filtering (for real-time applications) and smoothing (for a-posteriori analysis). Orbitography applications are detailed in [3]. These methods based on a differential correction process are strongly dependent on the measurement quality. Data pre-processing is necessary in order to get a reliable trajectory assessment. An additional issue occurs when several undiscriminated objects appear on the observations. This paper deals with such multiple trajectography problems.

---

[*] Ariane Group 78130 Les Mureaux, France
max.cerf@ariane.group

Problem description

A system of several independent objects is considered. Each object moves according to its own deterministic dynamical law, starting from an unknown initial state. An observation device takes photographs of the system at discrete dates. Each photograph exhibits a different number of objects and these objects are not discriminated (i.e. one does not know where each object appears on a given photograph). The trajectography goal is to reconstitute the greatest number of object trajectories from the series of photographs.

When a single object is observed, the trajectography aims at finding the trajectory parameters matching the measurements. For multiple objects, the available measurements cannot be associated a priori to the objects. An embedded assignment problem must be addressed at each observation date. The trajectography task (find the trajectory parameters) is complicated by an identification task (assign the measurements to the objects).

Figure 1 illustrates the problem for a set of 4 successive photographs showing respectively 3, 4, 4 and 2 objects. The object trajectories are modelled as straight lines.

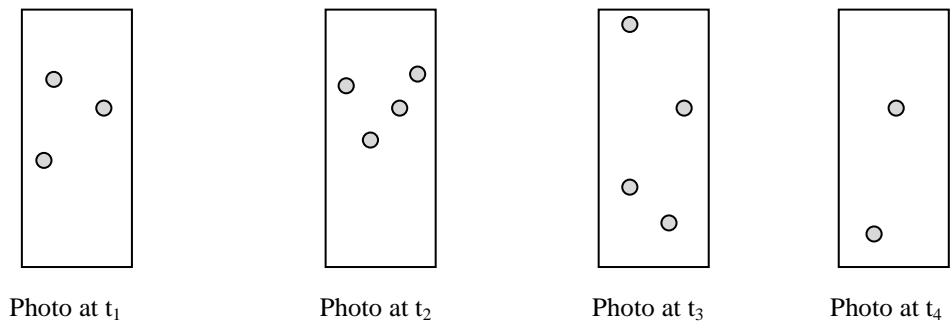

Figure 1 : Set of four photographs

A possible working hypothesis is to reconstitute two object trajectories from these observations. In that case the solution would have to discard one object from the photo 1 and two objects from the photos 2 and 3. This hypothesis could lead to the solution depicted on Figure 2.

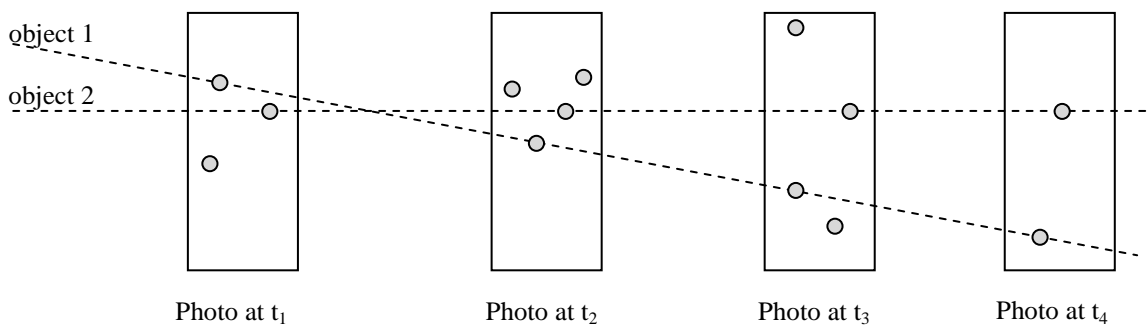

Figure 2 : Solution assuming two objects

Another possible working hypothesis would be to reconstitute four object trajectories, assuming that one object is missing on the photo 1 and two objects are missing on the photo 4.

The hypothesis could lead to the solution depicted on Figure 3. This solution identifies four objects matching the observations.

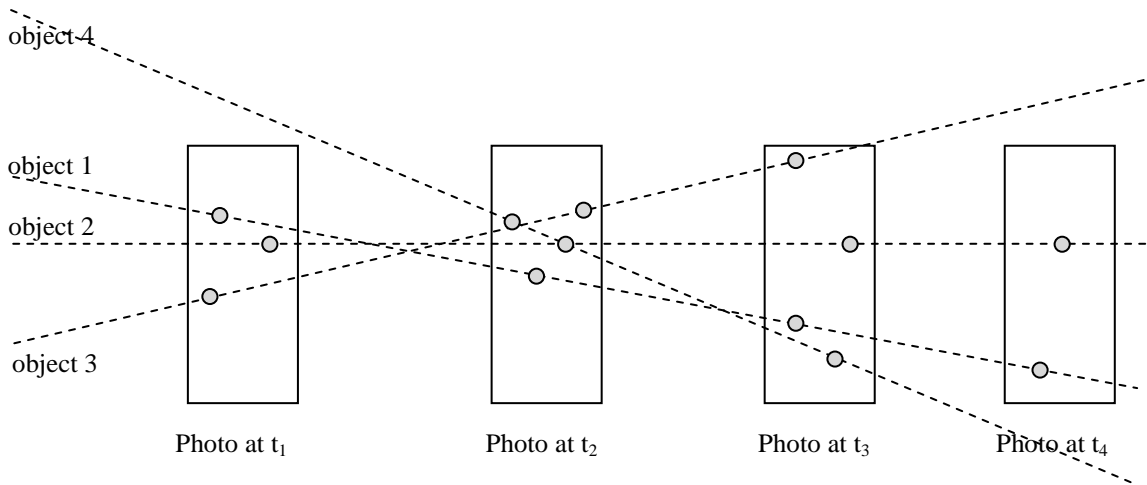

Figure 3: Solution assuming four objects

The solution assuming four objects makes the best possible use of the available observations, but such a solution may be difficult to achieve and less reliable for noisy observations. The conservative solution assuming only two objects is less sensitive to the observation defects. In this paper, the conservative approach is chosen based on the following working hypothesis.

Working hypothesis

The number of objects to identify is at most equal to the lowest number of objects observed per photograph.

This approach will yield a first trajectography result. The remaining observations may then be processed through a second trajectography run after discarding the observations selected by the first run. This iterative process proves more reliable than directly trying to identify as many objects as possible in one single shot.

The problem is mathematically formulated in §2. An algorithmic solution method is proposed in §3 and an application case is presented in §4.

## 2. Problem Formulation

This section presents the mathematical formulation of the multiple trajectography problem (MTP).

The number of objects to identify is denoted n. Each object is represented by its state vector x of dimension $n_x$. The evolution of the object number i ($1 \leq i \leq n$) is defined by the initial value problem

$$\begin{cases} \dot{x}_i = f_i(x_i, t) \\ x_i(t_0) = x_{i0} \end{cases}, \quad x_i \in R^{n_x}, \quad 1 \leq i \leq n \tag{1}$$

Each object evolves according to its own deterministic dynamics defined by the function $f_i$. The initial date $t_0$ is fixed prior to the first observation date. The initial states $(x_{i0})_{1 \leq i \leq n}$ are the problem unknowns. They are gathered in the vector $X_0$ of dimension $n_x \times n$.

$$X_0 = (x_{10}, \cdots, x_{n0}) \in R^{n_x \times n} \tag{2}$$

For a single object whose state vector is x at the date t, the observation device yields a measurement vector y of dimension $n_y$.

$$y = g(x, t) \in R^{n_y} \tag{3}$$

This equation assumes a perfect observation system. In a real application the measurement accuracy is limited. The uncertainties on the measurement vector y are introduced below.

A series of M observations is available at successive dates $t_1, \ldots, t_M$. Observing the object number i at the date $t_j$ would yield an expected measurement (also called "pseudo-measurement")

$$y_{ij} = g(x_i(t_j), t_j) \tag{4}$$

The set of the pseudo-measurements associated to the n objects at the date $t_j$ is gathered in the matrix $Y_j$ of dimension $n_y \times n$.

$$Y_j = (y_{1j}, \cdots, y_{nj}) \in R^{n_y \times n} \tag{5}$$

At the date $t_j$ the observation device takes a photograph, consisting in a set of simultaneous measurements made on the dynamical system. Every measurement of the photograph corresponds to an object.

In an ideal case, the number of measurements would be exactly n, and furthermore each object would be unambiguously identified. The trajectography could then be done independently for each object using its associated measurements.

In a real case, the number of measurements at the date $t_j$ is $m_j$. This number $m_j$ may change from one date to another. According to the working hypothesis stated in the introduction, the number n of targeted objects has been set to $n = \min(m_j)_{1 \leq j \leq M}$. This number $m_j$ is therefore either equal to or greater than the number n of targeted objects. The objects are not discriminated so that the measurements cannot be directly assigned to the targeted objects.

The $m_j$ measurements available at the date $t_j$ are denoted $(z_{kj})_{1 \leq k \leq m_j}$, where $z_{kj}$ is a vector of dimension $n_y$. They are gathered in the matrix $Z_j$ of dimension $n_y \times m_j$.

$$Z_j = \left( z_{1j}, \cdots, z_{m_j j} \right) \in \mathbf{R}^{n_y \times m_j} \tag{6}$$

The trajectography goal is to find the initial conditions $X_0$ and the resulting trajectories that match "as well as possible" the available measurements $Z_1, \ldots, Z_M$. The approach proposed in this paper consists in formulating an optimization problem with a "fitness function" defined as follows.

For given initial conditions $X_0$ the object trajectories are computed by integration of the dynamics equation Eq. (1). The set of pseudo-measurements $Y_j$ is then generated at each observation date $t_j$ using Eq. (4) and it must be compared to the set of available measurements $Z_j$.

Each measurement is a vector of dimension $n_y$. A norm denoted $\|.\|_j$ is defined on the measurement space at the date $t_j$ by

$$\|u\|_j = \sum_{k=1}^{n_y} \left( \frac{u_k}{\sigma_{kj}} \right)^2 \quad , \quad u \in \mathbf{R}^{n_y} \tag{7}$$

where   $u_k$ is the $k^{th}$ component of the measurement vector u

   $\sigma_{kj}$ is the measurement uncertainty on the component $u_k$ at the date $t_j$

This norm takes into account the accuracy of the measurement device depending on the date.

The cost $c_{kij}$ of assigning at the date $t_j$ the measurement $z_{kj}$ to the object number i is then assessed as

$$c_{kij} = \|z_{kj} - y_{ij}\|_j \tag{8}$$

A cost matrix $C_j$ is built gathering the assignment costs $c_{kij}$ of the measurements (k = 1 to $m_j$) to the objects (i = 1 to n) at the date $t_j$.

$$C_j = \begin{pmatrix} c_{11j} & \cdots & c_{1ij} & \cdots & c_{1nj} \\ \vdots & & \vdots & & \vdots \\ c_{k1j} & \cdots & c_{kij} & \cdots & c_{knj} \\ \vdots & & \vdots & & \vdots \\ c_{m_j 1 j} & \cdots & c_{m_j i j} & \cdots & c_{m_j n j} \end{pmatrix} \begin{matrix} 1 \\ \\ k \quad \text{measurement number} \\ \\ m_j \end{matrix} \tag{9}$$

with column labels: 1, i, n (object number)

A series of M cost matrices $C_1, \ldots, C_M$ associated to the measurement dates is thus defined. These matrices have the same number n of columns, but different numbers $m_j$ of rows (with $m_j \geq n$).

The assignment problem at the date $t_j$ consists in associating one measurement per object while minimizing the total assignment cost. The optimal assignment cost at $t_j$ is denoted $K_j$. The fitness function denoted F is the sum of the optimal assignment costs at the successive measurement dates.

$$F = \sum_{j=1}^{M} K_j \qquad (10)$$

The fitness function can be equivalently computed from the trajectory residuals. The residual $R_i$ of the trajectory number i sums the deviations of the trajectory with its assigned measurements using Eq. (8). It increases along the trajectory at each measurement date. The fitness function is thus equal to the sum of the residual final values.

$$F = \sum_{j=1}^{n} R_i \qquad (11)$$

The assessment process of the fitness function is depicted on Figure 4.

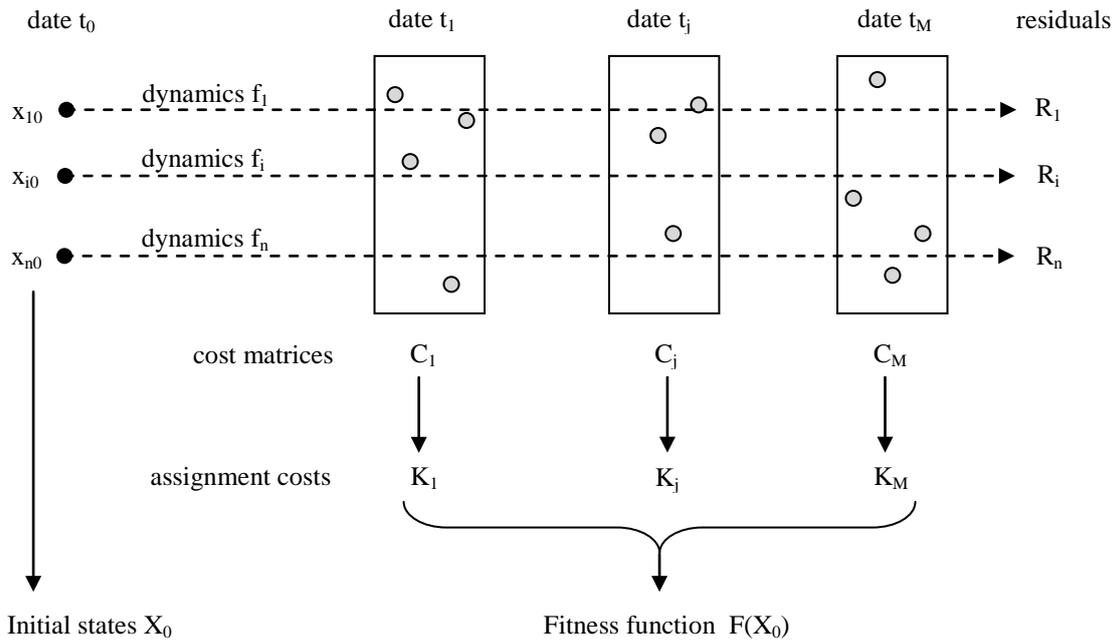

Figure 4 : Fitness function assessment

The fitness function depends on the object initial states $X_0$ which are the problem unknowns. The multiple trajectography optimization problem is formulated as

$$\min_{X_0} F(X_0) \qquad (12)$$

## 3. Solution algorithm

The multiple trajectory problem has been formulated as an optimization problem Eq. (12) with a fitness function to minimize. Assessing the fitness function requires solving a series of M embedded assignment problems with cost matrices $(C_j)_{1 \leq j \leq M}$ of variable sizes. Several solution methods may be envisioned, based either on general or on dedicated algorithms.

- A first method consists in formulating the assignment problem as an integer linear programming problem and solving it by a general linear programming algorithm (e.g. simplex, interior point). The constraint matrix of the assignment problem has the unimodularity property, which ensures finding directly an integer solution [4,5].
- A second method consists in formulating the assignment problem as a minimal cost flow problem in a bipartite graph. The problem can then be solved by the Busacker-Gowen algorithm or by another dedicated algorithm [6-8].
- A third method consists in applying the Hungarian algorithm of Kuhn-Munkres [9-11]. This dedicated algorithm amounts to solving the dual formulation of the assignment problem. The original method deals primarily with square cost matrices (same number of resources and tasks), but it can be easily adapted to non-square matrices. This algorithm is chosen for the MTP solution.

The fitness function F defined by Eq. (10) is not continuous since each assignment cost $K_j$ results from a combinatorial optimization problem. Gradient based algorithms are thus discarded and metaheuristic approaches are appropriate.

The Particle Swarm Optimization (PSO) method has proved successful on various engineering problems [12,13]. Several enhancements of the original method [14] have been proposed and tested, such as Repulsive Particle Swarm (RPS) [15] which is applied to the MTP. Each particle is defined as a vector of dimension $n_x \times n$ associated to a set of initial states $X_0$. The total number of particles is N. The particles move in the optimization search space (variables $X_0$). An iteration of the algorithm consists in modifying each particle velocity according to neighboring exchange rules and random perturbations, simulating the n object trajectories and assessing the resulting fitness value (Figure 4).

The algorithm proposed in [15] is applied with its standard settings. The main user options are the number of particles, the neighborhood size and topology (closest particles, ring, random, …) and the number of iterations. The velocity update formula is a weighted sum of individual, group and global components. A local search is performed by every particle along and opposite to its current velocity after moving. Optionally a part of the worst particles may also be randomly reinitialized.

## 4. Application case

The multiple trajectography algorithm is illustrated on a space application. Several objects evolving on near-geostationary Earth orbits are observed by a ground telescope. These objects may be for example debris stemming from a spent spacecraft or formation-flying satellites. The goal is to determine each object orbit from the available observations.

The algorithm is tested on an ideal case using synthetic measurements (opposite to "actual" measurements that would come from a real observation device). These synthetic measurements are obtained by simulating the object trajectories from given initial conditions and storing the resulting pseudo-measurements at discrete dates. Additional fictitious measurements are also manually added.

The optimization algorithm should theoretically be able to retrieve exactly the initial conditions and the trajectories of all objects.

A central gravity field without perturbations is considered for the object dynamics. The objects move on slightly eccentric orbits close to the geostationary ring (radius a = 42164 km) and on different orbital planes. The orbital plane is defined by the inclination I (angle with the Equator) and the right ascension of the ascending node Ω (RAAN) as depicted on Figure 5. The anomaly θ defines the initial angular position on the orbit from the ascending node [16].

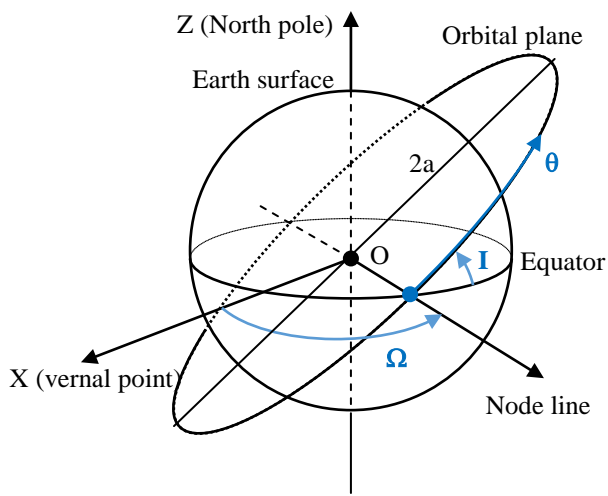

Figure 5 : Orbital parameters

a : semi-major axis
I : inclination
Ω : right ascension of the ascending node
   (RAAN)
θ : anomaly

A set of 10 objects is simulated on a duration of 3 days with the initial conditions given in Table 1.

|  | Semi-major axis a (km) | Eccentricity e (–) | Inclination I (deg) | RAAN Ω (deg) | Anomaly θ (deg) |
|---|---|---|---|---|---|
| Object 1 | 42164 − 100 | 0.02 | 0.5 | 0. | −40. |
| Object 2 | 42164 + 100 | 0.04 | 1.0 | 10. | −52. |
| Object 3 | 42164 − 50 | 0.06 | 0.6 | 20. | −64. |
| Object 4 | 42164 + 50 | 0.08 | 1.1 | 0. | −46. |
| Object 5 | 42164 − 20 | 0.03 | 0.7 | 10. | −58. |
| Object 6 | 42164 + 20 | 0.05 | 1.2 | 20. | −70. |
| Object 7 | 42164 − 80 | 0.02 | 0.8 | 0. | −42. |
| Object 8 | 42164 + 80 | 0.04 | 1.3 | 10. | −54. |
| Object 9 | 42164 − 40 | 0.06 | 0.9 | 20. | −66. |
| Object 10 | 42164 + 40 | 0.08 | 1.4 | 0. | −48. |

Table 1 : Initial conditions for the 10 objects

The observations are made by a ground telescope located at 0 deg longitude and 45 deg latitude. The telescope takes photographs of the sky every 30 minutes during 3 consecutive nights, each night extending over 5 hours. A set of 3×10 photographs is thus available for the trajectography.

The measurements derived from the photograph post-processing are the angles of elevation α (angle with the local horizontal plane) and azimuth ψ (angle with the North in the local horizontal plane, positive eastwards) defining the direction of the object seen from the telescope. The measurement geometry is depicted on Figure 6.

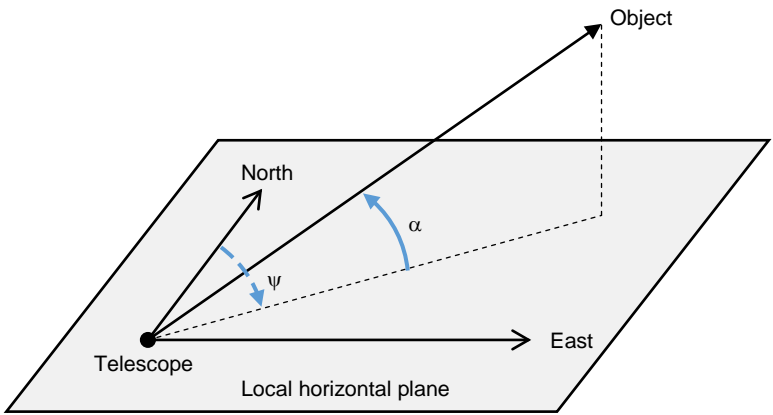

Figure 6 : Observation angles

α : elevation angle

ψ : azimuth angle

Figure 7 shows the pseudo-measurements obtained for each observation night. Each plot superposes the 10 consecutive photographs taken during one night. Some individual object trajectories may be visually guessed on these plots (on the right), but most objects are too close to be clearly discriminated.

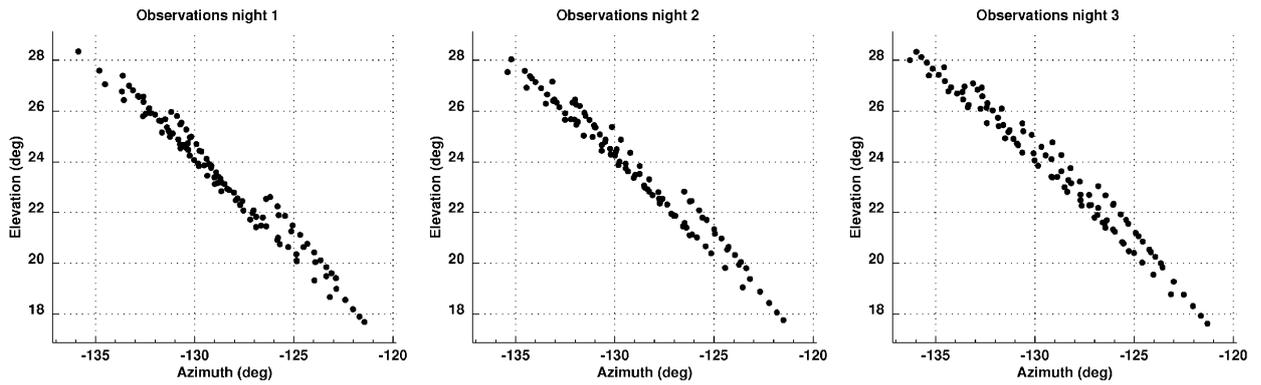

Figure 7 : Photographs taken during three consecutive nights

The elevation and azimuth angles associated to the 10 objects are plotted on Figure 8. The dashed lines represent the continuous evolution of the angles over the whole observation period (3 nights). The bold points mark the angles values picked up at the photograph dates. These values will be used as measurements for the algorithm test.

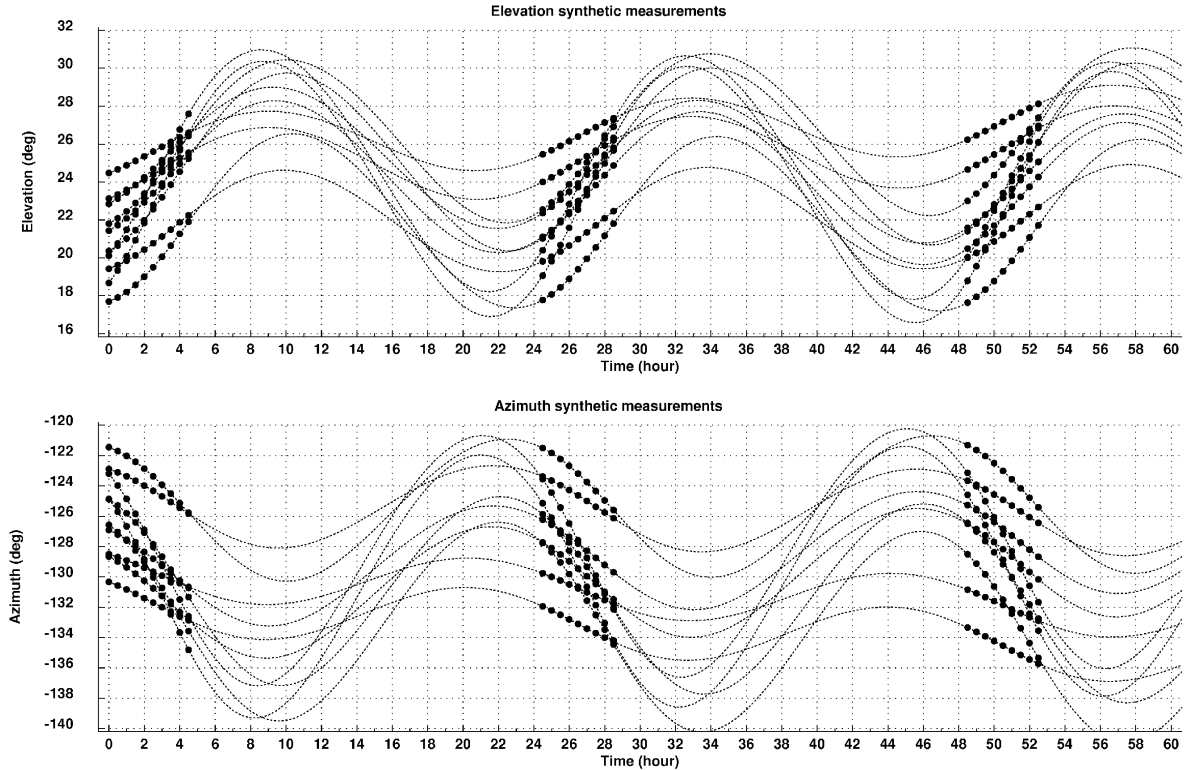

Figure 8 : Elevation and azimuth synthetic measurements

These initial simulations provide 30 matrices $(Y_j)_{1 \leq j \leq 30}$ (30 dates) of dimension 2×10 (2 angles measured, 10 objects). The measurement matrices $(Z_j)_{1 \leq j \leq 30}$ for the optimization test are obtained by adding one or two fictitious measurements at some dates, so that the number of objects appearing on each photograph ranges from 10 to 12. The uncertainty on both angles is set to 0.01 deg at all dates. The trajectography aims at finding 10 object trajectories matching the measurement matrices $(Z_j)_{1 \leq j \leq 30}$.

The particle swarm algorithm is run with the following settings : 100 particles, 400 iterations, 200000 functions calls, 20 closest neighbor topology, 2 local search steps at each iteration, 2 worst particles reset per iteration. The algorithm is parallelized is order to assess simultaneously all the particles.

The unknowns are the initial conditions for the 10 objects. The search bounds are :
- 42164 ± 200 km for the semi major axis (GEO altitude = 42164 km)
- [0 ; 0.1] for the eccentricity
- [0 deg ; 1.5 deg] for the inclination
- [-180 deg ; 180 deg] for the RAAN and the anomaly

The optimization takes about five minutes computation and it yields a set of 10 initial conditions presented in Table 2. For near equatorial orbits the RAAN is not well defined. The position is better defined by the total longitude $\lambda = \Omega + \theta$, which is given in Table 2. The object numbers in this table are relating to the variable order, they are thus not directly comparable to the numbers of Table 1.

|  | Semi-major axis a (km) | Eccentricity e (–) | Inclination I (deg) | Longitude $\lambda = \Omega + \theta$ (deg) |
|---|---|---|---|---|
| Object 1 | 42121 | 0.030 | 0.680 | -39.514 |
| Object 2 | 42199 | 0.054 | 1.078 | -46.769 |
| Object 3 | 42162 | 0.051 | 0.830 | -45.544 |
| Object 4 | 42168 | 0.055 | 0.945 | -44.533 |
| Object 5 | 42129 | 0.048 | 0.700 | -42.936 |
| Object 6 | 42192 | 0.037 | 0.919 | -49.164 |
| Object 7 | 42249 | 0.053 | 0.861 | -44.462 |
| Object 8 | 42152 | 0.035 | 0.650 | -41.394 |
| Object 9 | 42190 | 0.076 | 0.924 | -45.998 |
| Object 10 | 42115 | 0.043 | 1.009 | -49.606 |

Table 2 : Initial conditions found by optimization

Comparing these results to the "true" initial conditions of Table 1, which were used to generate the measurements, it is observed that the true initial conditions are not exactly retrieved. There are two explanations to this inaccuracy.

The first reason comes from the trajectography problem itself. Angles only orbit determination is indeed a challenging orbitography problem even for a single object [17]. All orbital parameters cannot be exactly retrieved when the observation time span is too small. For near-geostationary satellites observed from the Earth the problem is furthermore complicated by the distance assessment, which limits the achievable accuracy on the semi-major axis [18].

The second reason comes from the particle swarm algorithm. This algorithm is well suited to global optimization, but it is not adapted for converging accurately to a local minimum. Better results can be obtained by increasing the number of particles and the number of iterations, but at the expense of much larger computation times. In fact this is not useful for the multiple trajectography problem whose first goal consists in discriminating the objects, as analyzed here after.

Although the initial conditions are not exactly retrieved the measurement are correctly assigned. This is checked by the residual values $(R_i)_{1 \leq i \leq 10}$ of each individual trajectory. A low final residual indicates that the simulated trajectory is close to all assigned measurements. A low value of the fitness function Eq. (11) means that all final residuals are low.

Figure 9 shows the fitness function decrease during the iterations. A fast improvement is observed during the 10 first iterations, and a quite low value ($< 100.$) of the fitness function is finally achieved. Indeed 1 deg deviation on a single measurement, with an uncertainty set to 0.01 deg in Eq. (7), yields a contribution of 100 on the fitness function value.
The measurements residuals along the 10 optimized trajectories are also plotted on Figure 9. The residuals increase during the observation periods (plot as bold lines) and stay constant between two observation periods (plot as dashed lines). It can be observed that all trajectories contribute in a similar proportion to the global cost measured by the fitness function.

These results can be considered as satisfying although the initial trajectories are not perfectly retrieved. The primary goal of the multiple trajectory algorithm is indeed to correctly assign the measurements to the objects in order to discriminate them. A more accurate orbit determination using classical orbitography methods can then be achieved for each object considered separately.

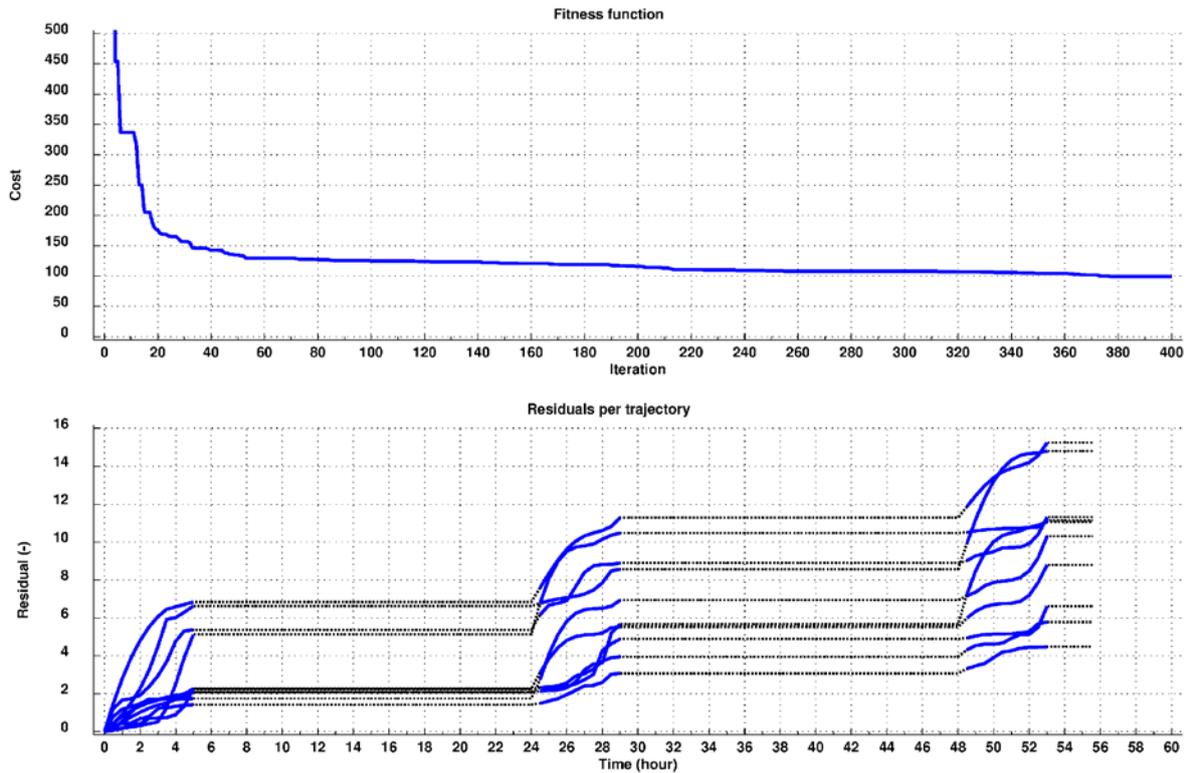

Figure 9 : Fitness function decrease and measurement residual per trajectory

## 5. Conclusion

Simultaneous trajectography of several objects raises assignment problems. The measurements available at each observation date must be correctly associated to the objects in order to yield reliable trajectory estimates. The problem is formulated as a global optimization problem solved by a particle swarm algorithm. Each particle represent a set of initial conditions used to simulate the object trajectories. An embedded assignment problem is solved at each observation date by the Hungarian method in order to associate the measurements to the simulations. The fitness function to minimize is the sum of all assignement costs. The method is illustrated on an orbitography example using optical observations of near geostationary satellites. It can be applied to various problems using for example radar or sonar measurements, with several objects to discriminate. More complex problems may be addressed by the same optimization approach, for example to estimate maneuvers or parameters of the dynamical model.

The approach proposed in this paper should be considered as a pre-processing stage of a multiple trajectography problem. It allows assigning correctly the measurements to the objects. The trajectography of each single object reduces then to a standard problem that can be addressed by usual filtering or smoothing methods in order to get a refined trajectory estimate.

# References

[1] Gelb A., *Applied Optimal Estimation*, Cambridge MA, MIT Press, 1989.
[2] Ljung L., *System identification: Theory for the user*, *Edition 2*, Prentive Hall, 1999.
[3] Vallado D.A., *Fundamentals of Astrodynamics and Applications*, Springer, 2007.
[4] Schrijver A., *Theory of Integer and Linear Programming*, Wiley and Sons, 1986.
[5] Nemhauser G., Wolsey L., *Integer and Combinatorial Optimization*, Wiley and Sons, 1999.
[6] Busacker R.G., Gowen P.J., *A procedure for determining minimal-cost network flow patterns*, ORO Technical Report 15, Operational Research Office, John Hopkins University, 1961.
[7] Orlin J.B., *A faster strongly polynomial minimum cost flow algorithm*, Proceedings of the $20^{th}$ ACM Symposium on the Theory of Computing, p. 377-387, 1988.
[8] Ravindra K.A., Thomas L.M., Orlin J.B., *Network Flows: Theory, Algorithms, and Applications*, Prentice-Hall, 1993.
[9] Kuhn H.W., *The Hungarian Method for the Assignment Problem*, Naval Research Logistics Quarterly, vol. 2, p. 83-97, 1955.
[10] Christofides N., *Graph Theory: An Algorithmic Approach*, Academic Press, 1975.
[11] Papadimitriou C., Steiglitz K., *Combinatorial Optimization*, Dover, 1998.
[12] Conway B., *Spacecraft Trajectory Optimization*, Cambridge University Press, 2010.
[13] Bonyadi M.R., Michalewicz Z., *Particle Swarm Optimization for Single Objective Continuous Space Problems: A Review*, Evolutionary Computation, vol. 25, issue 1, p. 1-54, 2017.
[14] Kennedy J., Eberhart R., *Particle swarm optimization*, IEEE International Conference on Neural Networks Proceedings, vol. 4, p. 1942-1948, 1995.
[15] Mishra S.K., *Repulsive Particle Swarm Method on Some Difficult Test Problems of Global Optimization*, MPRA Paper no. 1742, 2006.
[16] Chobotov, V., *Orbital mechanics*, $3^{rd}$ Edition, AIAA Education Series, 2002
[17] Escobal, P.R., *Methods of Orbit Determination*, $2^{nd}$ Edition, Wiley and Sons, 1976
[18] Fadrique, F.M., Maté A.A., Grau J.J., Sanchez J.F., Garcia L.A., *Comparison of angles only initial orbit determination algorithms for space debris cataloguing of Orbit Determination*, Journal of Aerospace Engineering, Vol IV, No 1, 2012

# Acronyms

| | |
|---|---|
| MTP | Multiple Trajectography Problem |
| PSO | Particle Swarm Optimization |
| RPS | Repulsive Particle Swarm |
| GEO | Geostationary Earth Orbit |
| RAAN | Right Ascension of Ascending Node |